\pdfoutput=1

\documentclass[10pt]{article}

\usepackage{geometry}
\usepackage{graphicx}
\usepackage{amsfonts}     
\usepackage[pdftex,bookmarks=false,colorlinks=true,pdfstartview=FitBV,linkcolor=blue,citecolor=blue,urlcolor=blue]{hyperref}        
\usepackage{amsmath} 
\usepackage{amssymb} 
\usepackage{amsthm}
\usepackage{subfigure}
\usepackage{cite}
\usepackage{tabularx}

\newcommand{\diff}[1]{\text{d}{#1}}

\newcolumntype{C}[1]{>{\centering\arraybackslash}p{#1}}

\renewcommand{\cite}[1]{[\citen{#1}]}

\DeclareMathOperator{\vecop}{vec}

\newtheorem{theorem}{Theorem}[section]

\newtheorem{problem}[theorem]{Problem}

\DeclareMathOperator{\kvec}{vec}

\title{C-DOC: Co-State Desensitized Optimal Control}

\author{
  Venkata Ramana Makkapati, Dipankar Maity, Mehregan Dor, Panagiotis Tsiotras\\
    School of Aerospace Engineering\\ 
    Georgia Institute of Technology\\ 
    Atlanta, GA 30332-0150 \\
  \texttt{\small \{mvramana, dmaity, mehregan.dor, tsiotras\}@gatech.edu} \\
}

\begin{document}

\maketitle

\begin{abstract}\label{sec:abstract}
In this paper, co-states are used to develop a framework that desensitizes the optimal cost. A general formulation for an optimal control problem with fixed final time is considered. The proposed scheme involves elevating the parameters of interest into states, and further augmenting the co-state equations of the optimal control problem to the dynamical model. A running cost that penalizes the co-states of the targeted parameters is then added to the original cost function. The solution obtained by minimizing the augmented cost yields a control which reduces the dispersion of the original cost with respect to parametric variations. The relationship between co-states and the cost-to-go function, for any given control law, is established substantiating the approach. Numerical examples and Monte-Carlo simulations that demonstrate the proposed scheme are discussed.
\end{abstract}


\section{Introduction} \label{sec:intro}

Obtaining robust solutions against parametric variations in optimal control problems is a requirement in various applications, particularly in the fields of aerospace and robotics. 
For many problems, it is essential that for a given performance criterion, one can also ensure minimal dispersion in the total cost, in spite of variations in the problem parameters.
A simple example would be the case of finding a minimum-time trajectory between two points (cities) on a map when a trajectory with travel time that is insensitive to traffic flow variations is desirable.
In this example, the trajectory designer is interested in finding a good trade-off between the minimum cost solution and a solution that is less sensitive to variations in the model parameters (the traffic flow in the previous example).
Solutions obtained by existing optimal control theory are model-dependent, and their resilience to variations in uncertain system parameters is not guaranteed.

A great deal of research has been performed in robust control that primarily focuses on the stabilization of systems under uncertainty and $H_\infty$ optimal control \cite{bhattacharyya1995robust,zhou1998essentials}. 
Previous attempts have primarily aimed at stability and performance criteria defined over an infinite horizon~\cite{petersen2014robust, wang1992robust}. 
Questions regarding the sensitivity of the trajectory (or the cost) explicitly and its implications on the performance have been largely overlooked.
In \cite{amato2001finite}, the authors have addressed the implications of parametric uncertainties over a finite horizon using linear matrix inequalities (LMIs), however, the problem does not address the sensitivity of the performance with respect to the parameters.
Traditionally, robust optimal control \cite{zhou1996robust, korkel2004numerical, Nagy2004411, lin1998optimalrobot} and feedback control synthesis \cite{HOROWITZ1963v} have been used to address the issue of parametric uncertainty, with an inherent trade-off between cost and robustness to
be decided. 
Indeed, the increased cost is incurred due to
additional control effort, in magnitude or over time. 
The main goal of \emph{desensitized optimal control}  (DOC) is to alleviate the additional effort induced onto the control feedback loop by, instead, picking a trajectory which is less sensitive to variations under parametric uncertainty.

Desensitization of a solution can include addressing the problem of minimizing variations:  
a) in the optimal trajectory; b) in the final state; or c) in the optimal cost under variation in model parameters, for a given optimal control problem. 
The former two cases were dealt with in our previous work \cite{ramana2018doc}.
Three different formulations that employ sensitivity functions were put forth which desensitize either entire trajectories or the state at a particular 
time (e.g., final time).

An approach to desensitize the cost for an optimal control problem with fixed final time is presented in this paper.
To this end, we recall that the co-states in an optimal control problem are a measure of the sensitivity of the value function with respect to the states along the optimal trajectory \cite{bryson, Peterson1974}. 
In this paper, we first prove that the co-states indeed capture sensitivity of the \emph{cost-to-go} function with respect to perturbations in the state given \emph{any} prescribed control law $u(t)$, not just the optimal one. 
Using this fact, a new approach to solve the DOC problem is presented.

Early work on trajectory sensitivity design include those of Winsor and Roy \cite{Winsor1970}, who developed a technique to design controllers that provide assurance for system performance under mathematical modeling inaccuracy. 
The feasibility of the technique was established with appropriate simulation results. 
However, their work has been restricted to linear systems.
Following that work, several approaches including sensitivity-reduction for linear regulators, using increased-order augmented system \cite{SUBBAYYAN_sensitivity}, 
modification of weighting matrix \cite{VERDE_sensitivity_reduction}, feedback \cite{Zames_optisensitivity, Chang_optidist_red}, and an augmented cost function \cite{Tang_sensitivity_1990, cheli2006observer}, were all thoroughly analyzed.
The approach of using an augmented cost function was further tested on the linear quadratic regulator (LQR) problem, which was later applied for active suspension control in passenger cars \cite{cheli2006observer}.

The work by Seywald et al. on desensitized optimal control makes use of sensitivity matrices to obtain an optimal open-loop trajectory that is insensitive to first-order parametric variations \cite{seywald1996desensitized, seywald2019}.
The proposed approach elevates the parameters of interest to system states, and defines a sensitivity matrix that provides the first-order variation in the states at time $t$, given the variation in the states at some time $t'$ ($t' \leq t$). 
An appropriate sensitivity cost is added to the existing cost function, and the dynamics of the sensitivity matrix is augmented in the system dynamics to solve the resulting optimal control problem.
The approach was later extended to optimal control problems with control constraints \cite{seywald2003desensitized}, and it was used to solve the Mars pinpoint landing problem \cite{shen2010desensitizing}. 
Some extensions to the landing problem include considering uncertainties in atmospheric density and aerodynamic characteristics \cite{Xu2015}, and using direct collocation and nonlinear programming \cite{li2011mars}.
However, the sensitivity matrix based approaches \cite{seywald1996desensitized, seywald2019} requires propagating the original states, the targeted parameters, and the elements in the sensitivity matrix, resulting to a total of $(n+\ell)^2 + n + \ell$ number of states.
An alternative approach was presented in \cite{ramana2018doc} where the dimensionality of the state-space for the augmented problem is reduced to $n + n\ell$, using traditional sensitivity functions.
Also, techniques based on sensitivity matrices have close connections with covariance trajectory-shaping, which was studied by Zimmer \cite{zimmer2005reducing} and Small \cite{small2010optimal}. 
A more detailed discussion on the existing methods for desensitized control can be found in \cite{weinmann2012uncertain}.

The contributions of this paper are summarized as follows: 
First, we mathematically prove (Theorem \ref{thm:main}) the fact that for any given control law, the co-states of an optimal control problem capture the first-order variations of the cost-to-go function given the variations in the system states. 
Second, we provide an approach to desensitize (reduce the variations of) the cost with respect to parametric variations by elevating the parameters to system states, and then penalizing the associated co-states of the optimal control problem.
Third, it is shown that the co-states and the sensitivity matrices (in Ref. \cite{seywald2019}) are related.
The approach is demonstrated on the Zermelo's navigation problem, and on a standard linear system thereby establishing its efficacy.

The rest of the paper is organized as follows. 
Section~\ref{sec:motivation} presents some preliminaries and formulates the problem addressed in this work. 
Section~\ref{sec:analysis} first presents a derivation of the relationship between the cost-to-go function and the co-states, and then provides an approach to solve the co-state based desensitized optimal control (C-DOC) problem. 
Section~\ref{sec:examples} demonstrates the proposed approach with two examples. 
It also contains important observations which point at some nuances of the approach. 
Section~\ref{sec:relation} discusses the relation between co-states and sensitivity matrices.
Section~\ref{sec:conclude} concludes the paper with some directions for future work. 


\section{Background and Preliminaries} \label{sec:motivation}

Consider a standard optimal control problem of the form 
\begin{align}
    \inf_u~~ &\mathcal{J}(u,p) \triangleq \phi(x(t_f),t_f) + \int_{t_0}^{t_f}L(x(t),u(t),t) \, \diff{t},
\label{eq:cost}
\end{align}
subject to
\begin{subequations}
\begin{align}
\dot{x} = f(x,p,u,t),& \quad x(t_0) = x_0, \label{eq:ini-dynamics},\\ 
\psi(x(t_f),&t_f) = 0, \label{eq:subjectto}
\end{align}
\end{subequations}
where $t \in [t_0,\,t_f]$ denotes time, with $t_0$ being the initial time and $t_f$ being the final time (both assumed to be fixed),  $p \in \mathcal{P} \subset \mathbb{R}^\ell$ are $\ell$ unknown, possibly time-varying, model parameters, $x(t) \in \mathbb{R}^n$ denotes the state, with $x_0$ being the fixed state at $t_0$. 
The control $u \in \mathcal{U}=\{ \text{Piecewise Continuous (PWC)},\,u(t) \in  U, ~ \, \forall\, t \in [t_0,t_f] \},$
with $U \subseteq \mathbb{R}^{m}$ being the set of allowable values of $u(t)$, $\phi:\mathbb{R}^{n}\times [t_0,t_f] \rightarrow \mathbb{R}$, the terminal cost function, and $L:\mathbb{R}^{n}\times\mathbb{R}^{m}\times [t_0,t_f] \rightarrow \mathbb{R}$, the running cost.
Finally, $\psi: \mathbb{R}^n \times [t_0,t_f] \rightarrow \mathbb{R}^k$ is a function representing $k$-number of constraint equations at the final time.
The above problem is to be solved by finding the optimal control $u^* \in \mathcal{U}$ that minimizes the cost function in (\ref{eq:cost}). 
The solution involves the optimal trajectory $x^*(t)$, $t \in [t_0,t_f]$, determined from $\dot{x}^*(t) = f(x^*(t), p, u^*(t), t)$ subject to $x^*(t_0) = x_0$, and the optimal cost $\mathcal{J}^* = \phi(x^*(t_f),t_f) + \int_{t_0}^{t_f}L(x^*(t),u^*(t),t) \, \diff{t}$.

The system dynamics represented by the function $f(x,p,u,t)$ contains the model parameters $p$ which are assumed to be constant.
It is understood that the optimal solution (constituting the cost, and the trajectory) is model-sensitive and, if changes in the parameters $p$ occur at any time $t \in [t_0,t_f]$, then the optimality of the obtained solution is not guaranteed.
Consequently, the optimal control problem has to be resolved
for each new value of the parameter vector.
If the optimal solution $u^*$ is used despite the parameter variations, one can expect a dispersion in the optimal trajectory (and/or cost $\mathcal{J}^*$).
With a motivation to minimize the dispersion of the final state $x^*(t_f)$ of the optimal solution, under parametric uncertainties, Seywald and Kumar constructed an augmented cost function using sensitivity matrices~\cite{seywald1996desensitized}.
The approach goes as follows.

First, the parameters of interest and the corresponding entries in the sensitivity matrix are elevated to states, and the augmented state $[ \tilde{x}^\top~(\vecop S)^\top] ^\top$,  where $\tilde{x} = [x~p]^\top$,
along with the corresponding dynamics and initial conditions are derived.
The sensitivity of the vector $\tilde{x}(t)$\footnote{From time to time we will denote $\tilde x(t)$ as $\tilde{x}(t|t_0,\tilde x_0)$ to explicitly represent the dependency on the initial conditions $\tilde{x}_0 = [x(t_0)^\top ~p(t_0)^\top]$} at time $t$ with respect to perturbations in the initial state vector $\tilde{x}(t_0) = \tilde{x}_0$ is denoted as $S(t|t_0,\tilde x_0)$
That is,
\begin{equation} \label{Sey:SF}
S(t|t_0,\tilde x_0) = \frac{\partial \tilde{x}(t|t_0,\tilde{x}_0)}{\partial \tilde{x}_0}.
\end{equation}
The dynamics of the  state $\tilde x$ can be written as
\begin{align}
&\dot{\tilde{x}} = \tilde{f}(\tilde{x},u,t) = [f^\top(x,p,u,t)~~0_{1 \times \ell}]^\top, \nonumber \\ 
&\tilde{x}(t_0) = \tilde{x}_0 = [x^\top_0~~p^\top_0]^\top, \label{dynamics_augstate}
\end{align}
and
\begin{align}
\dot{S}(t|t_0,\tilde x_0) = \frac{\partial \tilde{f}}{\partial \tilde{x}}S(t|t_0,\tilde x_0), \quad S(t_0|t_0,\tilde x_0) = I_{(n+\ell)}, \label{dynamics_seywald}
\end{align}
where  $p_0$ is the nominal value of the parameter vector, and where $S(t|t_0,\tilde x_0)$ represents the sensitivity of the vector $\tilde{x}(t)$ at time $t$ with respect to perturbations in the initial state vector $\tilde{x}(t_0)$.

The augmented cost function, given in (\ref{eq:seywald_aug}) below, is then minimized to obtain an optimal solution with the final state being ``desensitized" with respect to the parameter variations
\begin{align}     \label{eq:seywald_aug}
\mathcal{J}_s(u,p) = \mathcal{J}(u,p) + \int_{t_0}^{t_f} \| \vecop \left( S(t_f|t_0,\tilde x_0) S(t|t_0,\tilde x_0)^{-1}\right) \|_{Q(t)}^{2} \Big) \diff{t}, 
\end{align}
with $Q(t) \geq 0,~\text{for all}~t \geq t_0$.
Note that the sensitivity matrix of Seywald in (\ref{Sey:SF}) has the form of a state transition matrix and its properties are exploited to construct the sensitivity of the final state with respect to the variations in the state at time $t\in [0, t_f]$, which is then plugged into the running cost (\ref{eq:seywald_aug}).
This is elaborated upon in Ref. \cite{seywald1996desensitized}. The approach achieves the desensitization of the final state.

\subsection{Problem Formulation}

In this work, we are interested in desensitizing the cost itself. 
By denoting
\begin{align}
    \mathcal J(u,p)=\int_{t_0}^tL(x(s),u(s),s)\diff{t}+C(\tilde x(t),u,t),\\
    C(\tilde x(t),u,t)=\int_{t}^{t_f}L(x(s),u(s),s)\diff{t}+\phi(x(t_f),t_f),
\end{align}
we immediately notice that the parametric variation at time $t$, affects the total cost $\mathcal{J}(u,p)$ only through the cost-to-go $C(\tilde x(t),u,t)$.
Thus, the sensitivity of the total cost for a parametric variation at time $t$ from its nominal value $p_0$ can be captured through the term
\begin{align}\label{eq:d-cost}
   S_C(x(t),p_0,u,t)= \frac{\partial C}{\partial p}(\tilde x(t),u,t)\Big|_{p=p_0}.
\end{align}
There are several ways to capture the effect of the parametric variations on the cost, one of which is to consider the following sensitivity cost
\begin{align}
   \mathcal{J}_c(u,p_0)= \int_{t_0}^{t_f}\|S_C(x(t),p_0,u,t)\|^2_{Q(t)}\diff{t},
\end{align}
for some $Q(t)\ge 0$, for all $t\ge t_0$.

There are three major formulations relevant to the problem of cost-based desensitization, which are as follows.

\begin{problem} \label{prob:1}
Solve
\begin{subequations}
\begin{align}
    \inf_{u\in \mathcal{U}}~~ & \mathcal{J}(u,p_0)\\
    \mathrm{subject~to}~~&\mathcal{J}_c(u,p_0)\le D.
\end{align}
\end{subequations}
\end{problem}

\vspace{2ex}

Let us denote the solution of Problem \ref{prob:1} to be the ``cost-desensitization'' function $\mathsf{J}(D)$ which represents the optimal cost given a bound on the sensitivity metric.
A similar problem is to consider minimizing the sensitivity of the cost for a given bound on the performance index, as presented below.

\begin{problem} \label{prob:2}
Solve
\begin{subequations}
\begin{align}
    \inf_{u\in \mathcal{U}}~~ & \mathcal{J}_c(u,p_0)\\
    \mathrm{subject~to}~~&\mathcal{J}(u,p_0)\le J.
\end{align}
\end{subequations}
\end{problem}

\vspace{2ex}

Let us denote the solution of Problem \ref{prob:2} to be the ``desensitization-cost'' function $\mathsf{D}(J)$. 
Finding analytical or numerical solutions to $\mathsf{J}(D)$ or $\mathsf{D}(J)$ are challenging.
However, $\mathsf{J}(D)$ or $\mathsf{D}(J)$ can be constructed by solving the following family optimization problems for all $\alpha\in [0,\infty)$.

\begin{problem} \label{prob:3}
Solve
\begin{align}
    \inf_{u\in \mathcal{U} }\mathcal{J}(u,p_0)+\alpha\mathcal{J}_c(u,p_0)
\end{align}
\end{problem}

\vspace{2ex}

By observing that the scalar $\alpha$ can be absorbed into the matrix $Q(t)$, we will rewrite the objetive function in Problem \ref{prob:3} as
\begin{align*}
    \mathcal{J}_s(u) =\mathcal{J}(u,p_0)+\mathcal{J}_c(u,p_0).
\end{align*}
When the sensitivity cost has zero weight ($Q(t)\equiv 0$), we solve problem \eqref{eq:cost} and retrieve $\limsup_{D\to \infty}\mathsf{J}(D)$, and as we increase the weight on the sensitivity cost (through $Q(t)$), we arrive at an optimal control whose performance is more insensitive to the variations in the parameters. 
In the limit when $Q(t)\to \infty$ for all $t$, we retrieve $\limsup_{J\to \infty}\mathsf{D}(J)$.
In this work, we will focus on minimizing $\mathcal{J}_s(u)$.
Detailed analysis of $\mathsf{J}(D)$ and $\mathsf{D}(J)$ will appear elsewhere.

The new optimization problem  we are interested in solving is
\begin{subequations}
\begin{align}
    \inf_u ~~&\mathcal{J}_s(u) \\
    \text{subject to}~&\dot{x}=f(x,p_0,u,t),~~x(t_0)=x_0,\\
    &\psi(x(t_f),t_f)=0. \label{eq:ter_const}
\end{align}
\end{subequations}

The following section presents a formal proof for the fact that the co-states capture the sensitivity of the cost-to-go function for any given control input $\bar{u}(t)$, that satisfies the terminal constraint (\ref{eq:ter_const}) with nominal value of the parameter $p_0$.
The result would allow us to penalize a weighted norm of the co-states, with their dynamics obtained from the adjoint equations, that desensitizes the cost function with respect to the variations in the targeted parameters.


\section{Co-states and Desensitized Optimal Control} \label{sec:analysis}
In this section we characterize the cost-sensitivity $S_C(x(t),p_0,u,t)$ in terms of the co-state process associated with the optimal control problem given by \eqref{eq:cost}-\eqref{eq:subjectto}.
The following theorem shows that the sensitivity of the cost-to-go function with respect to the state at time $t$ can be represented by a co-state process $\lambda$ with certain boundary conditions at the final time $t_f$.
\vspace{1ex}
\begin{theorem} \label{thm:main}
	Consider the dynamical system $\dot x=f(x,u,t)$, evolving under a given control law $\bar{u} \in \bar{\mathcal{U}}\subseteq \mathcal{U}$, where
	\begin{align*}
	\bar{\mathcal{U}} =  \Big\{ \bar{u}:[t_0, t_f] \rightarrow \mathbb{R}^{m}  \; \text{is}  \; \text{PWC}, \; \bar{u}(t) \in U, \text{ such that } \psi(x(t_f),t_f)=0, \; x(t_f) =x_0 + \int_{t_0}^{t_f}f(x(t),\bar{u}(t),t) \, \diff{t}  \Big\}.
	\end{align*}
Then, for a cost-to-go function (associated with the cost functional \eqref{eq:cost}) with $x = x(t)$
	\begin{align}
	C(x,\bar{u},t) = \phi(x(t_f),t_f) + \int_{t}^{t_f} L(x(\tau),\bar{u}(\tau),\tau) \, \diff{\tau},
	\label{eq:cost-to-go}
	\end{align}
	under the control $\bar{u} \in \bar{\mathcal{U}}$, the sensitivity of the cost-to-go function with respect to the state $x$ at time $t$ is, 
	\begin{equation}
	\lambda^{\top}(t) = \dfrac{\partial C}{\partial x}(x(t),\bar{u},t),
	\end{equation}
	which obeys the dynamics
	\begin{equation}
	\dot{\lambda}^{\top}(t) = - \frac{\partial H}{\partial x}(x(t),\bar{u},\lambda(t),t), \label{eq:thm_co-state_dyn}
	\end{equation}
	where
	\begin{align}
	H(x,u,\lambda,t) = L(x,u,t) + \lambda^{\top}f(x,u,t).
	\label{eq:Hamiltonian}
	\end{align}
	Furthermore, the terminal condition for (\ref{eq:thm_co-state_dyn}) is given by
	\begin{align}
	\lambda(t_f) = \frac{\partial \phi}{\partial x}(x(t_f),t_f). \label{trans_cond1}
	\end{align}
\end{theorem}

\vspace{1ex}

\begin{proof}
    The proof is presented in Appendix \ref{Ap:proof}.
\end{proof}

\vspace{1ex}

It is interesting to note that the theorem holds not only for the optimal control (a result that follows directly from the maximum principle \cite{pontryagin2018mathematical}), but for any control law that is piecewise-continuous and ensures that the terminal constraint is met. The C-DOC problem can now be fully formulated using this result.


For the C-DOC problem the augmented state is $\tilde{x}=[x^\top,p^\top]$ with dynamics given in \eqref{dynamics_augstate}.
The Hamiltonian, defined in Theorem \ref{thm:main}, for this system, can be written as
\begin{align}
H(\tilde{x}, u, \lambda, \mu, t) &= L(x,u,t) + \lambda^\top \dot{x} + \mu^\top \dot{p} \nonumber \\
&= L(x,u,t) + \lambda^\top f(x,p,u,t),
\end{align}
where $\lambda$ and $\mu$ are the co-states corresponding to state $x$ and vector of parameters defined by $p$, respectively. The corresponding adjoint equations are given by
\begin{align}
\dot{\lambda}^\top &= -\frac{\partial H}{\partial x}(\tilde x, u,\lambda,\mu, t) = -\lambda^\top \frac{\partial f}{\partial x}(x,p,u,t)-\frac{\partial L}{\partial x}(x,u,t), \label{eq:dyn_lambda} \\
\dot{\mu}^\top &= -\frac{\partial H}{\partial p}(\tilde x, u,\lambda,\mu, t) = -\lambda^\top \frac{\partial f}{\partial p}(x,p,u,t). \label{eq:dyn_mu}
\end{align}
Since the co-states represent the sensitivity of the cost-to-go function for a given control input $u(t)$ (Theorem \ref{thm:main}), they can be expressed as 
\begin{align}
\lambda(t)^\top &= \frac{\partial C}{\partial x}(\tilde x(t),u,t), \label{eq:lambdat}\\
\mu(t)^\top &= \frac{\partial C}{\partial p}(\tilde x(t),u,t), \label{eq:mut}
\end{align}
for a given control  $u\in \bar{\mathcal{U}}$, this results in the trajectory $x(t)$ for $t_0 \leq t \leq t_f$, where 
\begin{align*}
C(\tilde x,u,t) = \phi(x(t_f),t_f) + \int_{t}^{t_f}L(x(\tau),u(\tau),\tau) \, \diff{\tau}.
\end{align*}
Note that $p$ is an augmented state in the given problem and affects the cost $\mathcal{J}$ through the state $x$, whose dynamics is a function of $p$.
Since we have used $\dot p=0$ and $p(t_0)=p_0$, we have ensured that $p(t)=p_0$.
Thus, by comparing equations (\ref{eq:d-cost}) and (\ref{eq:mut}), we obtain $\mu(t) = S_C(x(t),p_0,u,t)$.
Therefore, weighting the co-state in the existing cost function will ensure that the solution of the augmented problem minimizes the sensitivity of the cost $\mathcal{J}$ with respect to parametric variations. This results in an updated optimal control problem with an augmented cost, accounting for the sensitivity component, given by 
\begin{align}
{\mathcal{J}}_s(u) &= \phi(x(t_f),t_f) + \int_{t_0}^{t_f}  \left[ L(x(t),u(t),t)  + \mu^\top(t) Q(t) \mu(t) \right] \, \diff{t}. \label{eq:desen_cost}
\end{align}
Minimizing the cost (\ref{eq:desen_cost}) subject to the dynamics (\ref{dynamics_augstate}), terminal constraint (\ref{eq:subjectto}), and the transversality conditions (\ref{trans_cond1}) with
\begin{align}
\mu(t_f) = 0,
\end{align}
yields a desensitized optimal control problem for the original problem. Here, $Q(t)\in \mathbb{R}^{\ell\times \ell}$ is a user-defined positive semi-definite weighting function and is generally of the form
\begin{align}
Q(t) \equiv \text{diag}(\alpha_1(t),\ldots, \alpha_{\ell}(t)).
\end{align}
  This co-state based approach requires formulating $2(n+\ell)$ number of states, as compared to the higher $2(n+\ell)^2 + n+ \ell$ states in \cite{seywald1996desensitized}, employing sensitivity matrices for an optimal control problem. 
The resulting problem \eqref{eq:desen_cost} is typically solved by the off-the-shelf existing solvers. 


\section{Numerical Examples} \label{sec:examples}

In many applications, the resulting trajectories should be insensitive with respect to perturbations and/or uncertainties within the model parameters at specified times along the trajectory. The following section presents some numerical examples that will aid in understanding the implementation of this technique and will elucidate its subtleties. The simulations are obtained using GPOPS-II \cite{Patterson2014}. 

\subsection{Zermelo's Navigation Problem}

Consider a typical Zermelo's problem \cite{seywald1996desensitized} with currents parallel to the shore ($x_1$) as a function of $x_2$ such that
\begin{align}
v_{\text{current}} = p x_2,
\end{align}
where $p$ is a parameter, which is uncertain, and its nominal value is $p_0 = 10$. The problem has to be desensitized with respect to this parameter. The dynamics can be written as
\begin{align} 
\dot{x}_1 &= \cos(u) + p x_2, \\
\dot{x}_2 &= \sin(u),
\end{align}
subject to the conditions
\begin{align*}
x_1(0)=0, \quad x_2(0)=0, \\
x_2(t_f)=0, \quad  t_f=1.
\end{align*}
For the case where the state $x_1(t_f)$ has to be maximized, the cost function can be expressed in the Mayer form as
\begin{align}
\min_{u} \mathcal{J}(u) = -x_1(t_f),
\end{align} 
where $u \in \mathcal{U}=\left\lbrace \text{C}^0, u(t) \in  [0,2\pi), \; \forall\ t \in \left[0,1\right]\right\rbrace$ and is the control input. 
Here, optimization of the trajectory is done with respect to parameter $p$, which is not precisely known. 
In order to facilitate desensitization of the cost (in this case the final state $x_1(t_f)$) 
with respect to variations in the parameter $p$ we first consider the augmented dynamics
\begin{align}
\dot{\tilde{x}} = \begin{bmatrix}
\dot{x}_1\\\dot{x}_2\\\dot{p}
\end{bmatrix}
=\begin{bmatrix}
\cos(u) + p x_2\\\sin(u)\\0
\end{bmatrix},
\end{align}
with boundary conditions 
\begin{align}
\tilde{x}(0) = \begin{bmatrix}
x_1(0)\\x_2(0)\\p(0)
\end{bmatrix} = \begin{bmatrix}
0\\0\\p_0
\end{bmatrix}, \quad x_2(t_f) = 0.
\end{align}
\begin{figure}[htb!]
	\centering
	\subfigure[Optimal trajectories]{\includegraphics[width=0.4\textwidth]{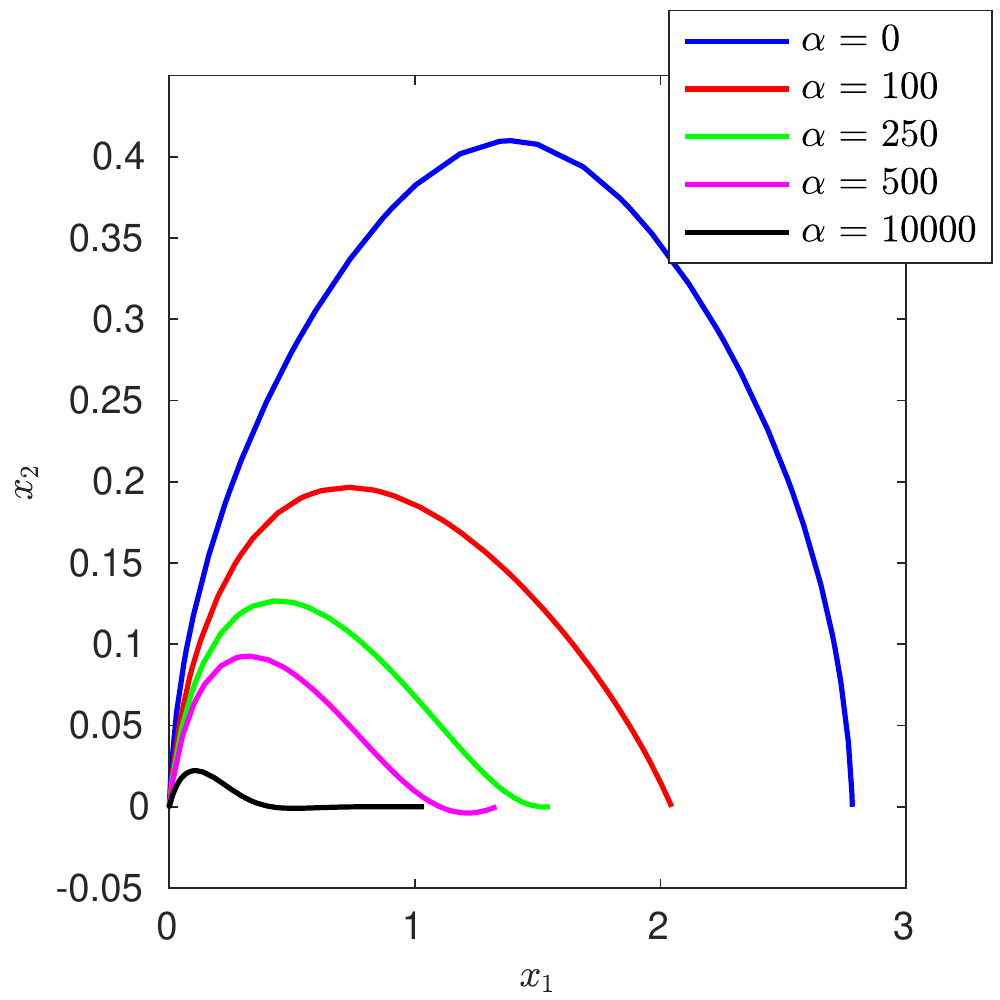}} 
	\subfigure[$\mu^2(t)$ - a measure of sensitivity]{\includegraphics[width=0.4\textwidth]{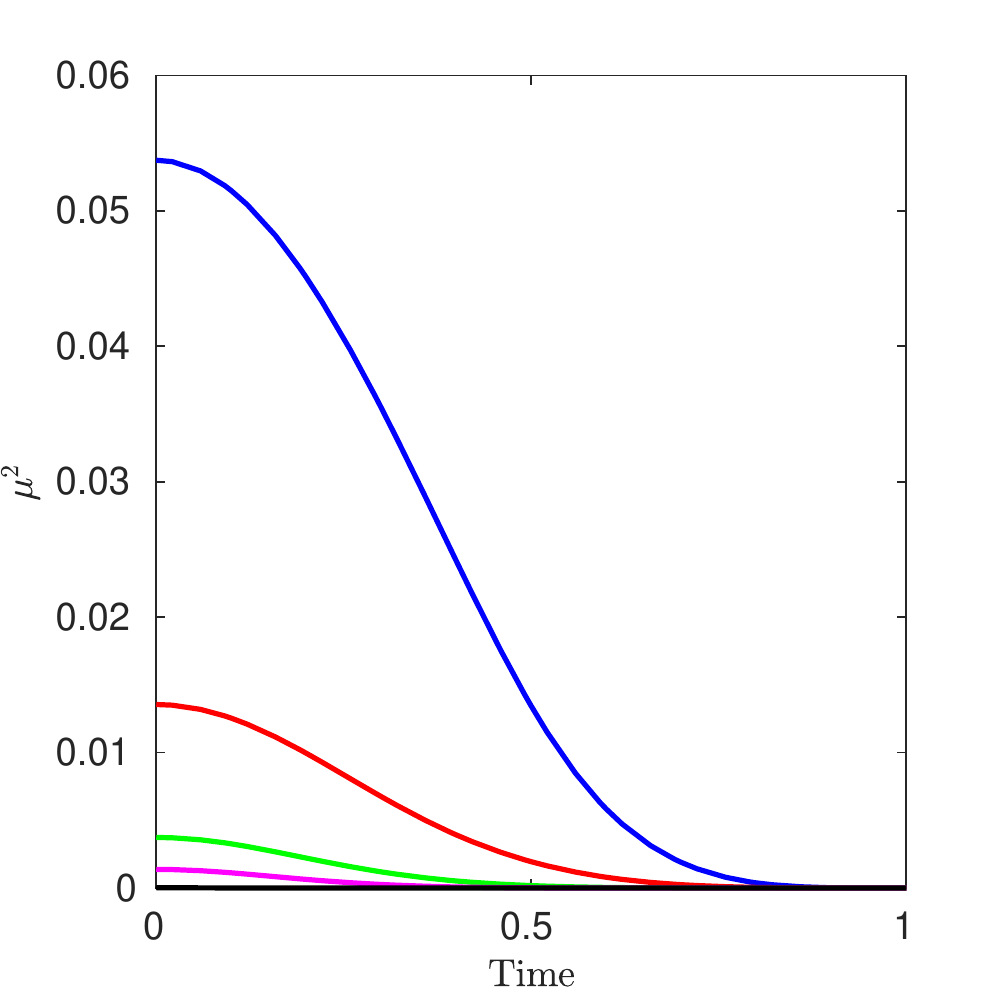}} \\
	\subfigure[Monte-Carlo simulations]{\includegraphics[width=0.4\textwidth]{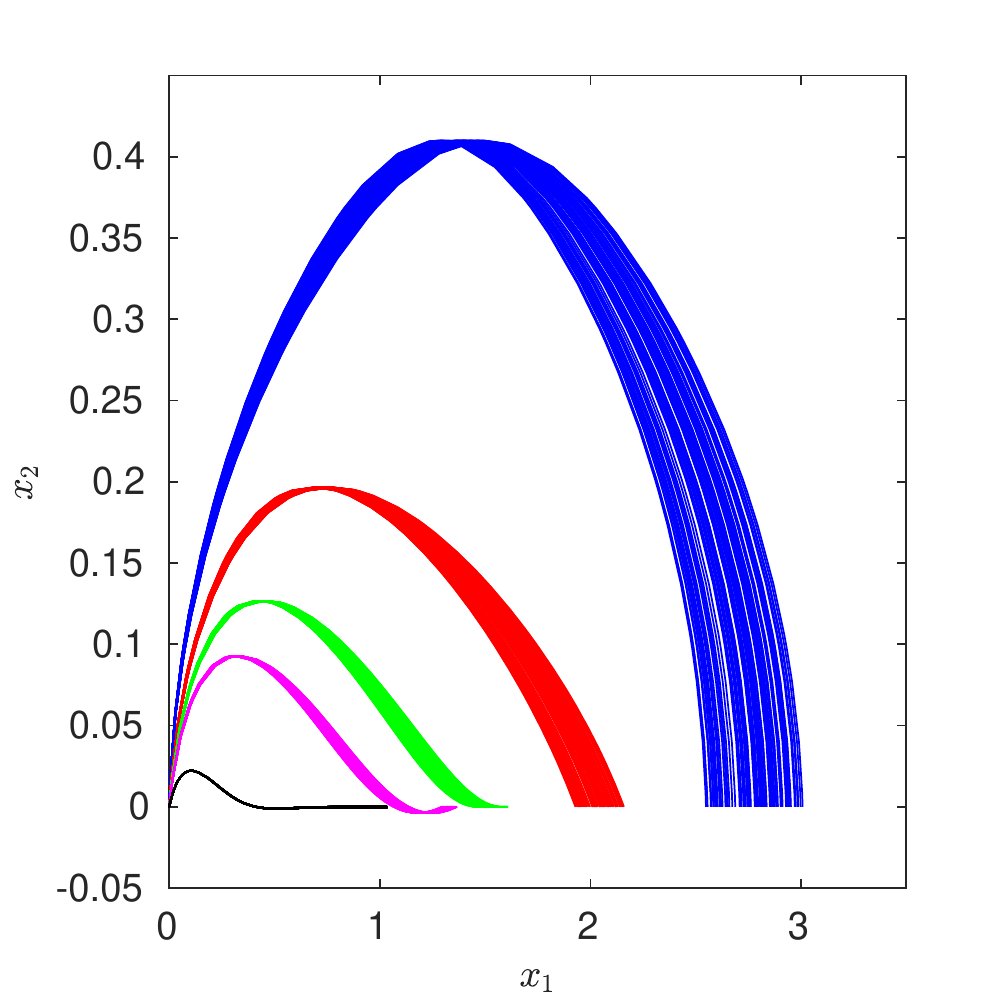}}
	\subfigure[Optimal Control ($u^*$)]{\includegraphics[width=0.4\textwidth]{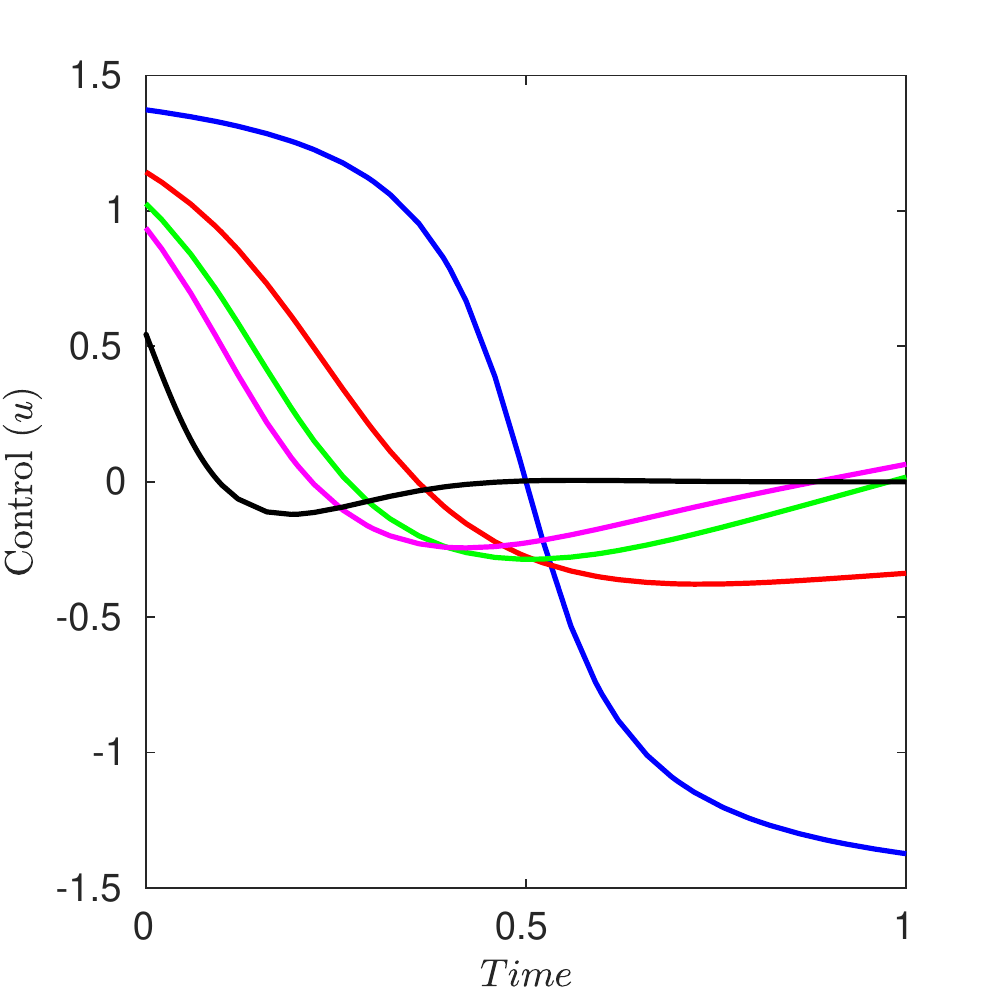}}  
	\caption{Results obtained for the Zermelo's path optimization problem with different levels of desensitization}
	\label{fig:Zermelo}
\end{figure}
The Hamiltonian is given by
\begin{align}
H(\tilde{x},u,\lambda, \mu) &= \lambda_1 (\cos(u)+px_2) + \lambda_2 \sin(u),
\end{align}
and the adjoint equations are:
\begin{align}
\dot{\lambda}_1 &= -\frac{\partial H}{\partial x_1} = 0, \\
\dot{\lambda}_2 &= -\frac{\partial H}{\partial x_2} = -\lambda_1 p, \\
\dot{\mu} &= -\frac{\partial H}{\partial p} = -\lambda_1 x_2.
\end{align}
The cost for the desensitized optimal control problem is
\begin{align}
\min_{u} {\mathcal{J}}_s(u) = -x_1(t_f) + \int_{0}^{t_f} \alpha \mu^2(t) \, \diff{t}.
\end{align}

The weight parameter $\alpha$ is chosen between 0 and 10,000. Figure \ref{fig:Zermelo}(a) shows the optimal paths obtained for the five different values of $\alpha$ in the selected range.
As the value of $\alpha$ increases, the optimal solution moves closer to the shore, minimizing the effect of uncertainty in the current, which has an effect proportional to the distance $x_2$. 
For $\alpha \rightarrow \infty$, the optimal path would be moving along the shore, i.e., along the $x_1$ axis.
The results obtained from the Monte-Carlo simulations are shown in Figure \ref{fig:Zermelo}(c) which confirm the expected desensitization of the cost $x_1(t_f)$, further substantiating the claims.
For the Monte-Carlo simulations, a time-constant parametric variation is enforced, and 100 different values between $[0.9p_0,1.1p_0]$ are randomly chosen to be $p$.
The trends of the integrand in the running cost $\mu^2(t)$ (Figure \ref{fig:Zermelo}(b)), show that its value is steadily decreasing and is almost zero for $\alpha = 10,000$. 

\subsection{Linear Systems}

Consider an optimal control problem of minimizing a quadratic cost 
\begin{align}
\mathcal{J}(u) = \int_{0}^{t_f} {\frac{1}{2}(x^\top R_1 x + u^\top R_2 u)} \, \diff{t},
\end{align}
given the $n$-dimensional linear dynamics with parameter vector $p$
\begin{align}
\dot{x} &= A(p)x + B(p)u,\\
\dot{p} &= 0,
\end{align}
with initial conditions
\begin{align}
x(0) = x_0, \quad p(0) = p_0,
\end{align}
where $x \in \mathbb{R}^n$, $u \in \mathbb{R}^m$, $p \in \mathbb{R}^\ell$, $A:\mathbb{R}^\ell \rightarrow \mathbb{R}^{n \times n}$, $B:\mathbb{R}^\ell \rightarrow \mathbb{R}^{n \times m}$, $R_1 \geq 0$, $R_2 > 0$, and $t_f$ is fixed. The goal is to desensitize the cost with respect to the parameter $p$. Following the steps to construct the cost term for desensitization, the Hamiltonian is given by
\begin{align}
H &= \frac{1}{2}(x^\top R_1 x + u^\top R_2 u) + \lambda^\top \dot{x} + \mu^\top \dot{p}, \nonumber \\
&= \frac{1}{2}(x^\top R_1 x + u^\top R_2 u) + \lambda^\top (A(p)x + B(p)u).
\end{align}
The adjoint equations are
\begin{align}
\dot{\lambda}^\top = -\frac{\partial H}{\partial x} &= -x^\top R_1 - \lambda^\top A(p),\\
\dot{\mu}^\top = -\frac{\partial H}{\partial p} &= - (x^\top \otimes \lambda^\top) \frac{\partial}{\partial p} \kvec(A(p)) - (u^\top \otimes \lambda^\top) \frac{\partial}{\partial p} \kvec(B(p)).
\end{align}
where $\lambda$ and $\mu$ are the co-states of $x$ and $p$, respectively. 
Since the cost has to be desensitized with respect to $p$, 
the augmented cost that has to be minimized for the C-DOC problem is given by
\begin{align}
{\mathcal{J}}_s(u) = \int_{0}^{t_f} {\frac{1}{2}(x^\top R_1 x + u^\top R_2 u + \mu^\top Q \mu)} \, \diff{t}.
\end{align}

\begin{figure}[htb!]
	\centering
	\subfigure[Optimal trajectories]{\includegraphics[width=0.4\textwidth]{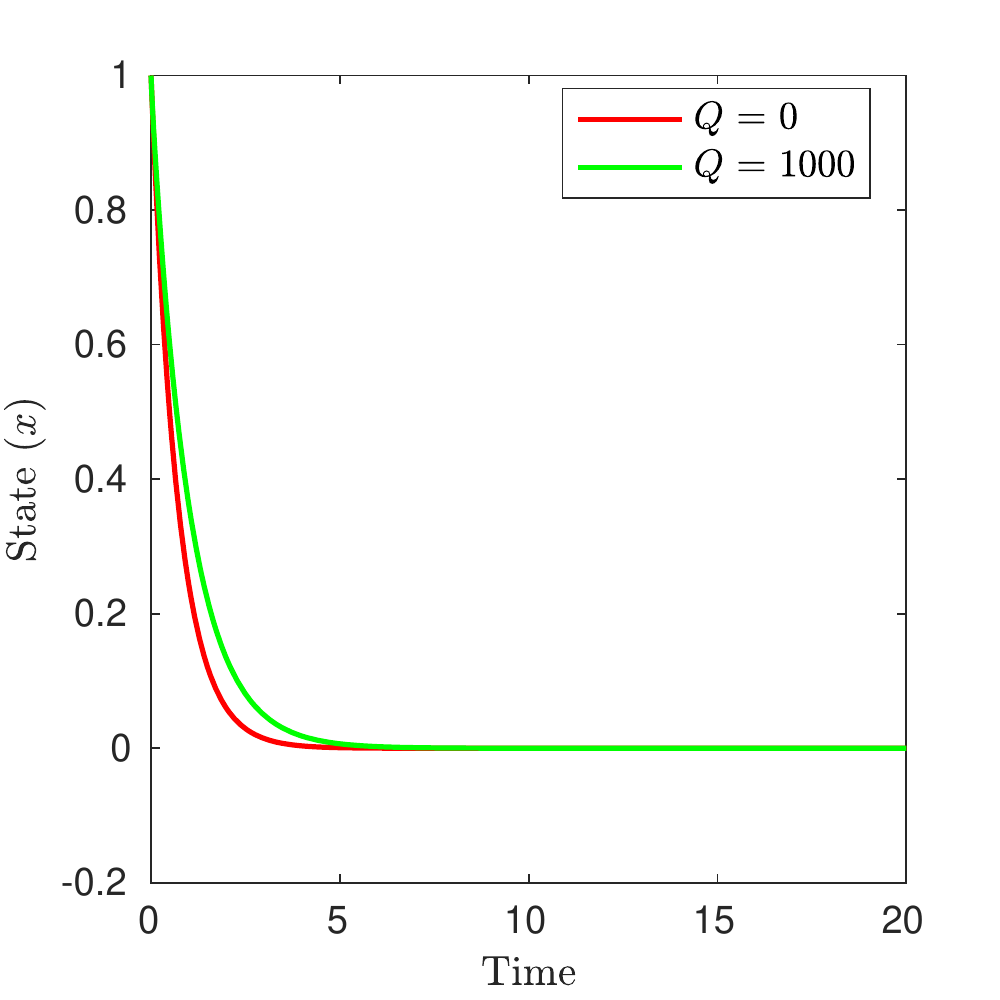}} 
	\subfigure[$\mu^2(t)$ - a measure of sensitivity]{\includegraphics[width=0.4\textwidth]{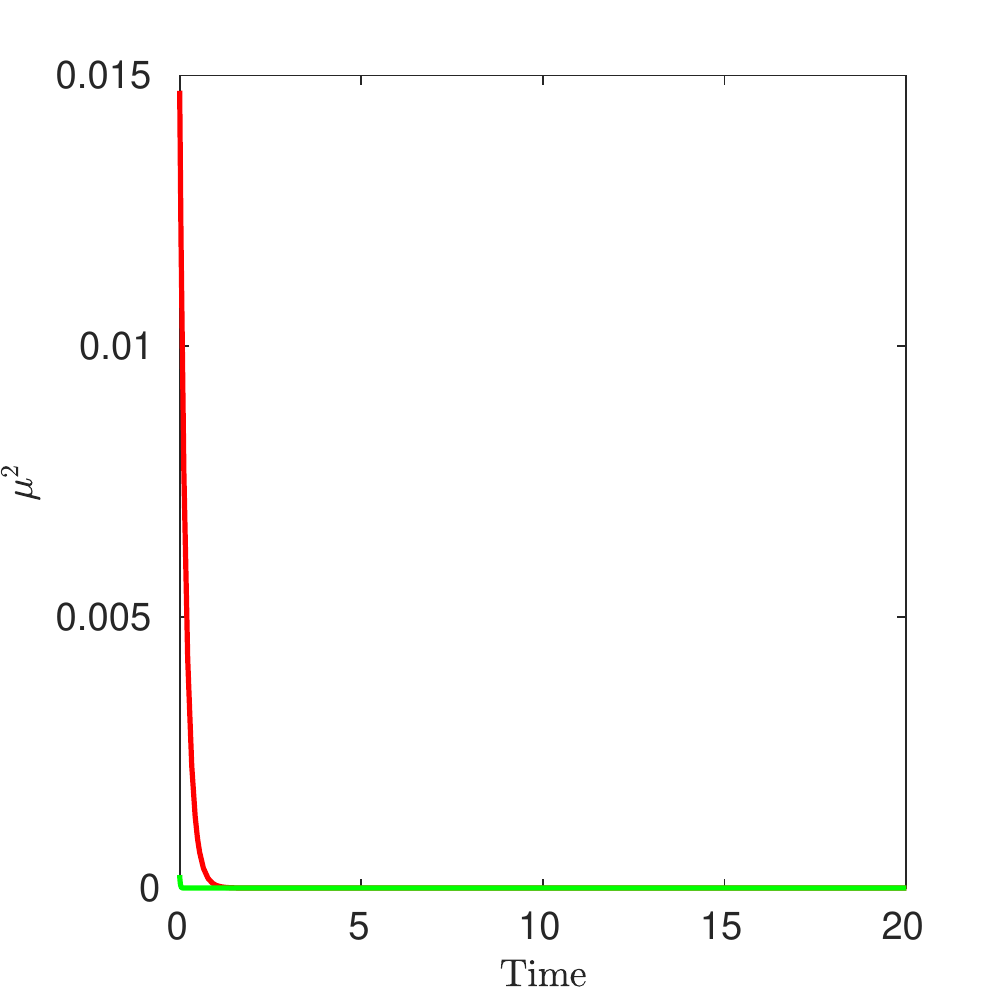}} \\
	\subfigure[Monte-Carlo simulations]{\includegraphics[width=0.4\textwidth]{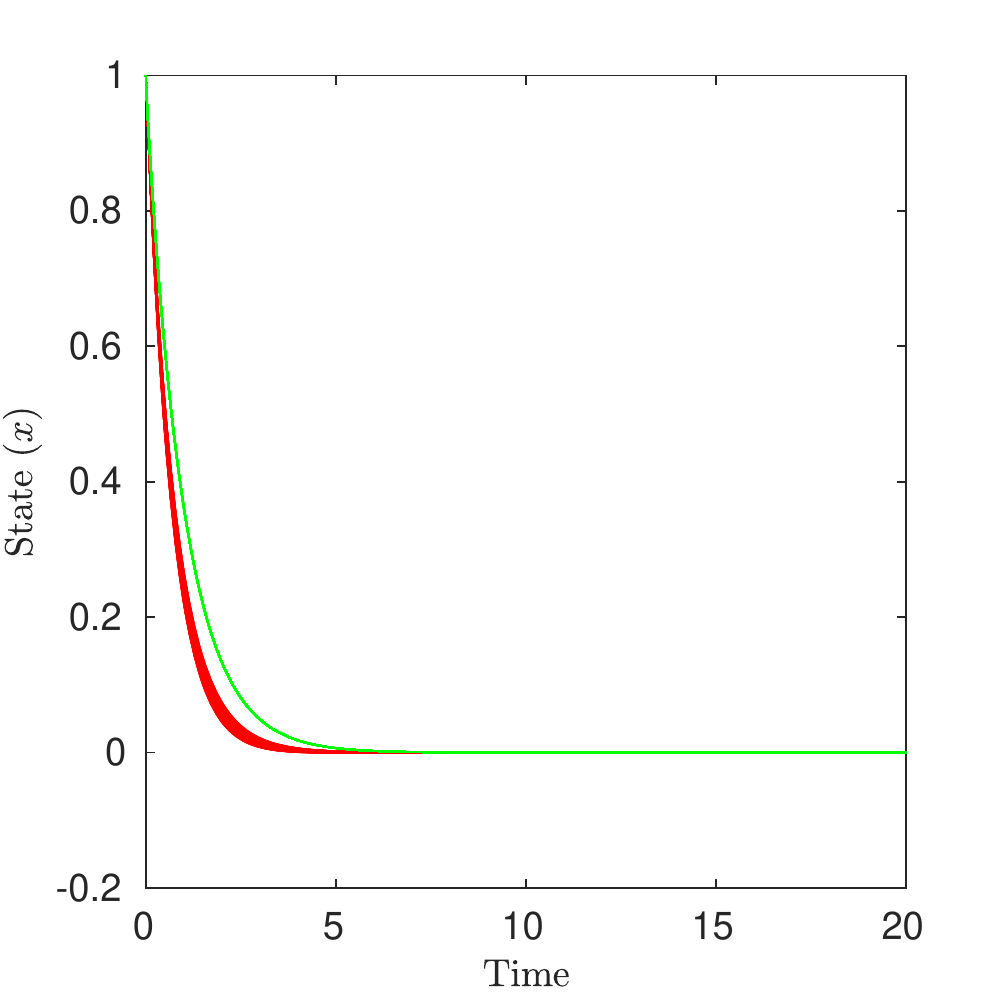}}
	\subfigure[Optimal Control ($u^*$)]{\includegraphics[width=0.4\textwidth]{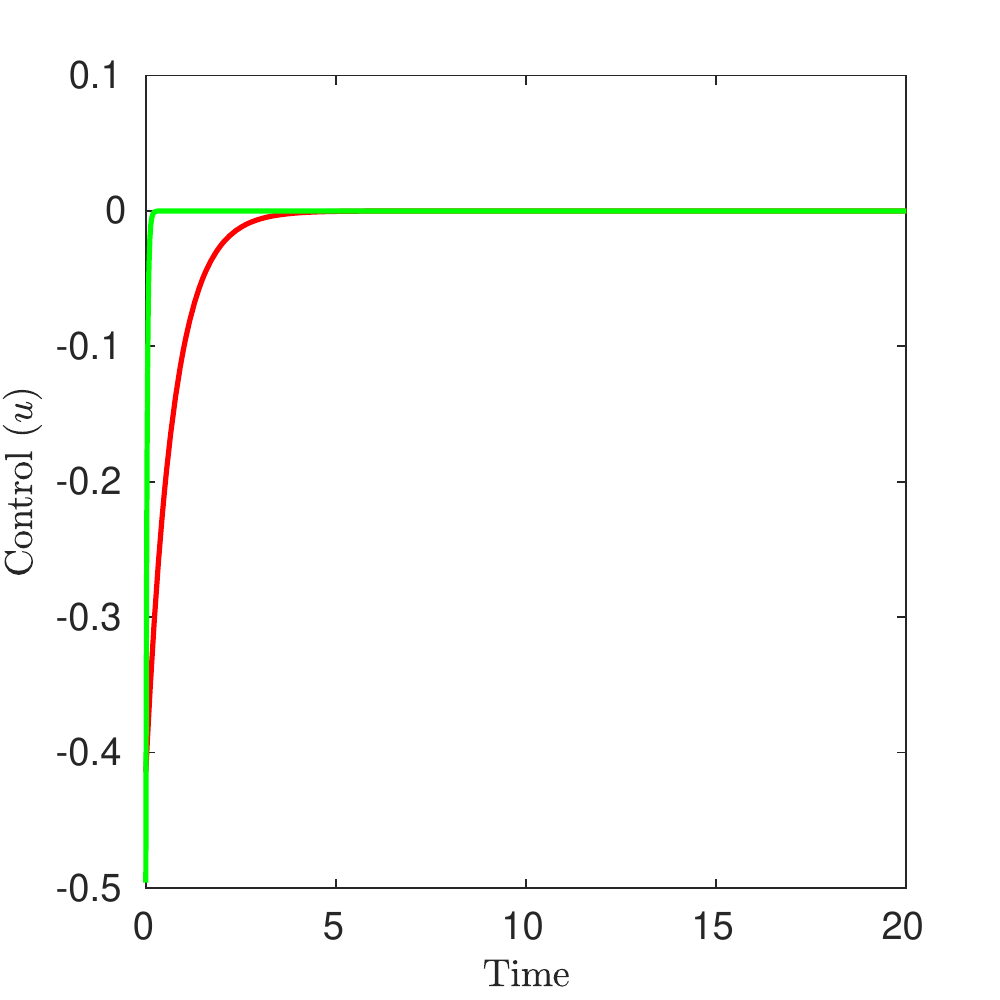}}  
	\caption{Results obtained for an LQR problem with $B$ matrix being uncertain}
	\label{fig:LS_b}
\end{figure}

To demonstrate the results, we consider a one-dimensional
linear system with the dynamics 
$\dot{x} = ax+bu$ with initial condition $x(0) = 1$, and let $R_1 = R_2 = 2$, $t_f = 20$.
We first analyze the case where $b$ is the uncertain parameter with its nominal value as $b_0 = 1$, and $a = -1$.
The solutions obtained for $Q = 0$ and $1,000$ are shown in Figure~\ref{fig:LS_b}.
Note that the sensitivity measure ($\mu^2(t)$) in Figure~\ref{fig:LS_b}(b) is lower for the desensitized solution.
Since $b$ is the source of uncertainty that perturbs the trajectory (and eventually the cost), by 
introducing desensitization ($Q = 1,000$), it can be observed from Figure~\ref{fig:LS_b}(d) that the control goes to zero earlier compared to the non-desensitized solution.
By making the control zero, the source of uncertainty is removed from the system.
The results obtained from the Monte-Carlo simulations with $ b \in [0.8b_0, 1.2b_0]$ are shown in Figure~\ref{fig:LS_b}(c), which suggests that the variation in the cost for the desensitized solution is significantly lower. 

\begin{figure}[htb!]
	\centering
	\subfigure[Optimal trajectories]{\includegraphics[width=0.4\textwidth]{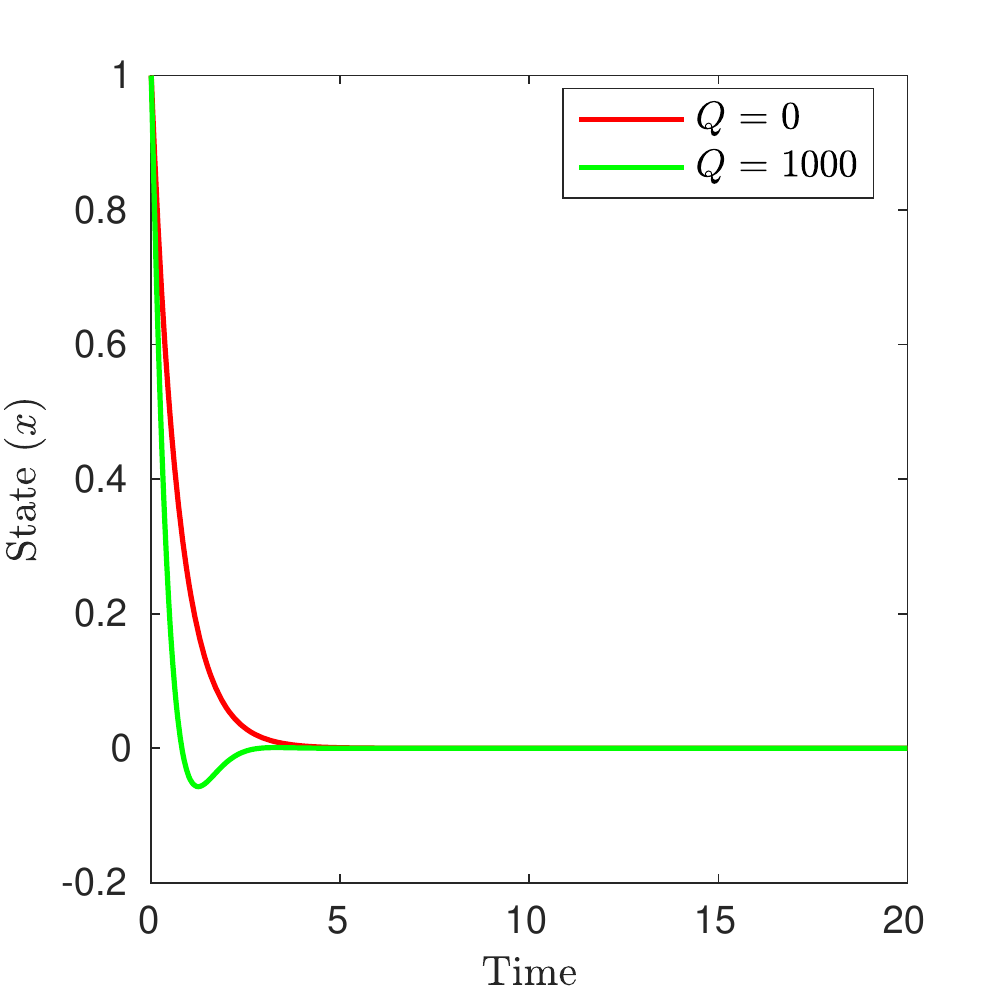}} 
	\subfigure[Optimal Control ($u^*$)]{\includegraphics[width=0.4\textwidth]{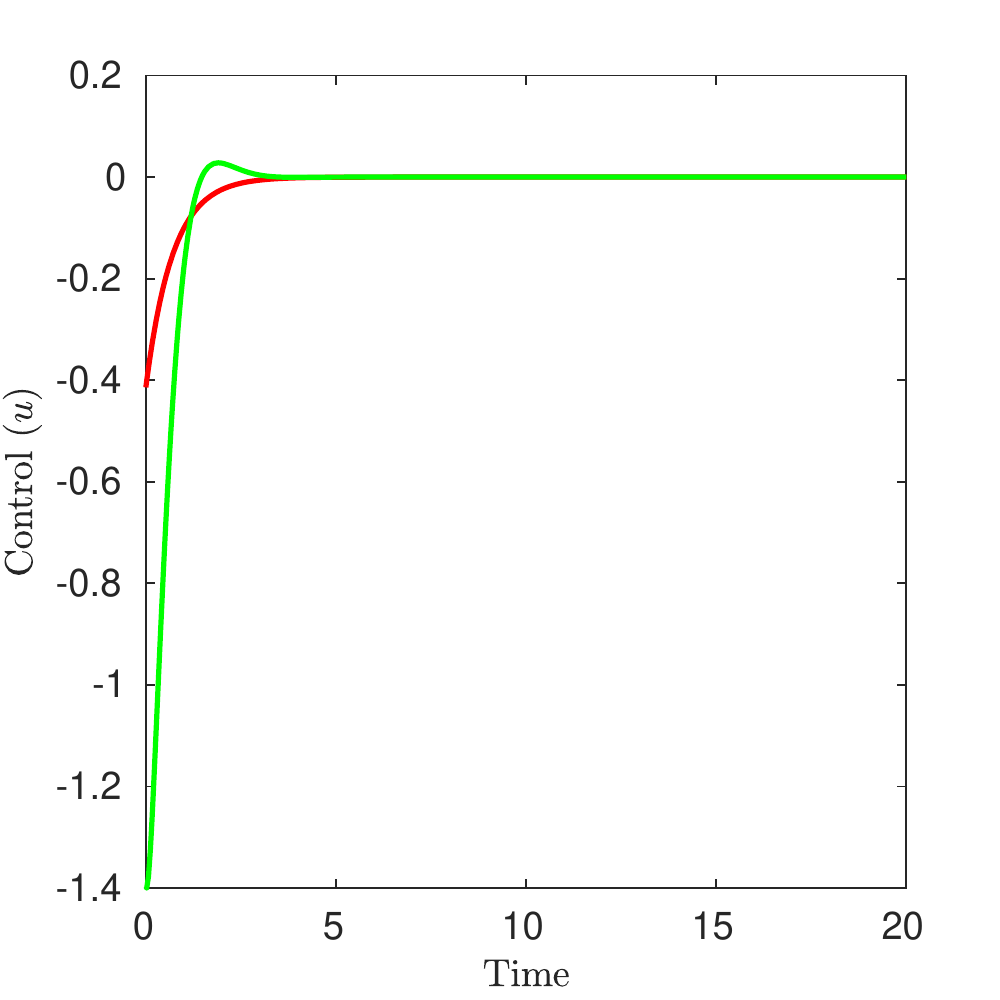}} \\
	\subfigure[Monte-Carlo simulations]{\includegraphics[width=0.4\textwidth]{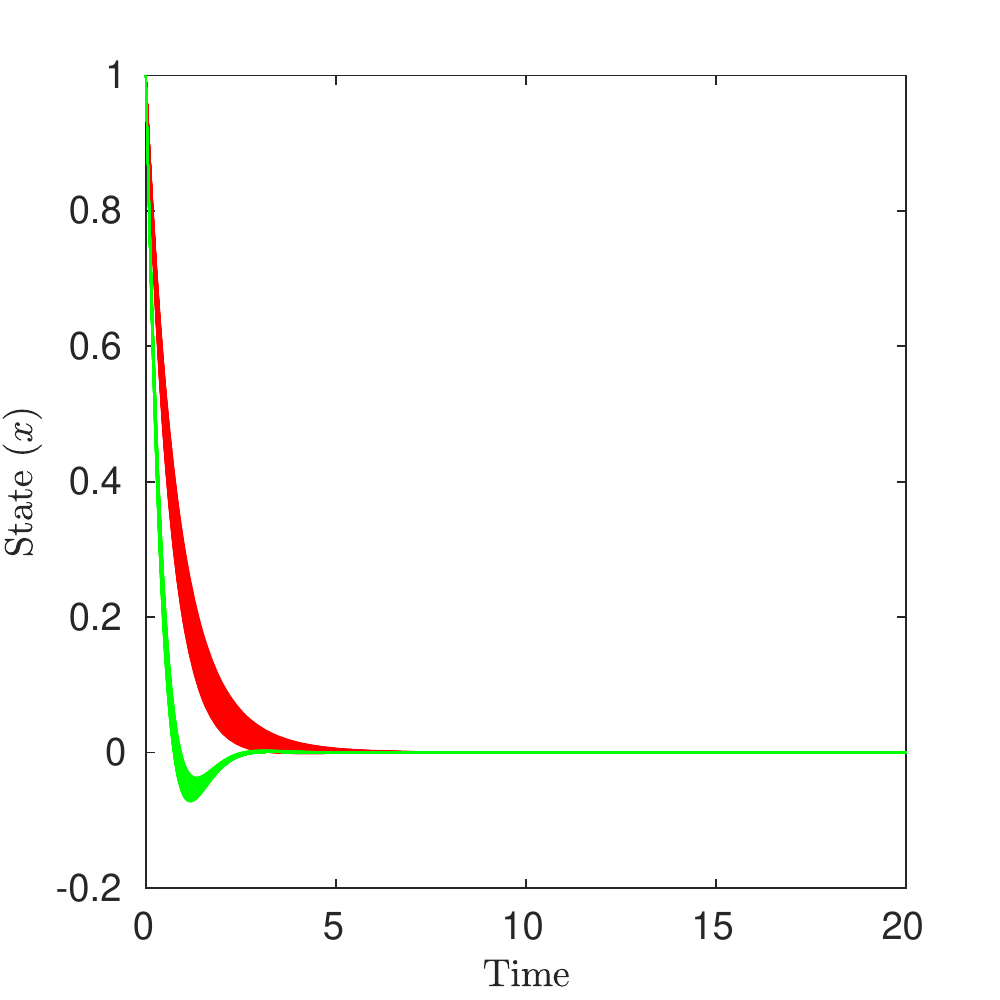}}
	\subfigure[Cost ($J$)]{\includegraphics[width=0.4\textwidth]{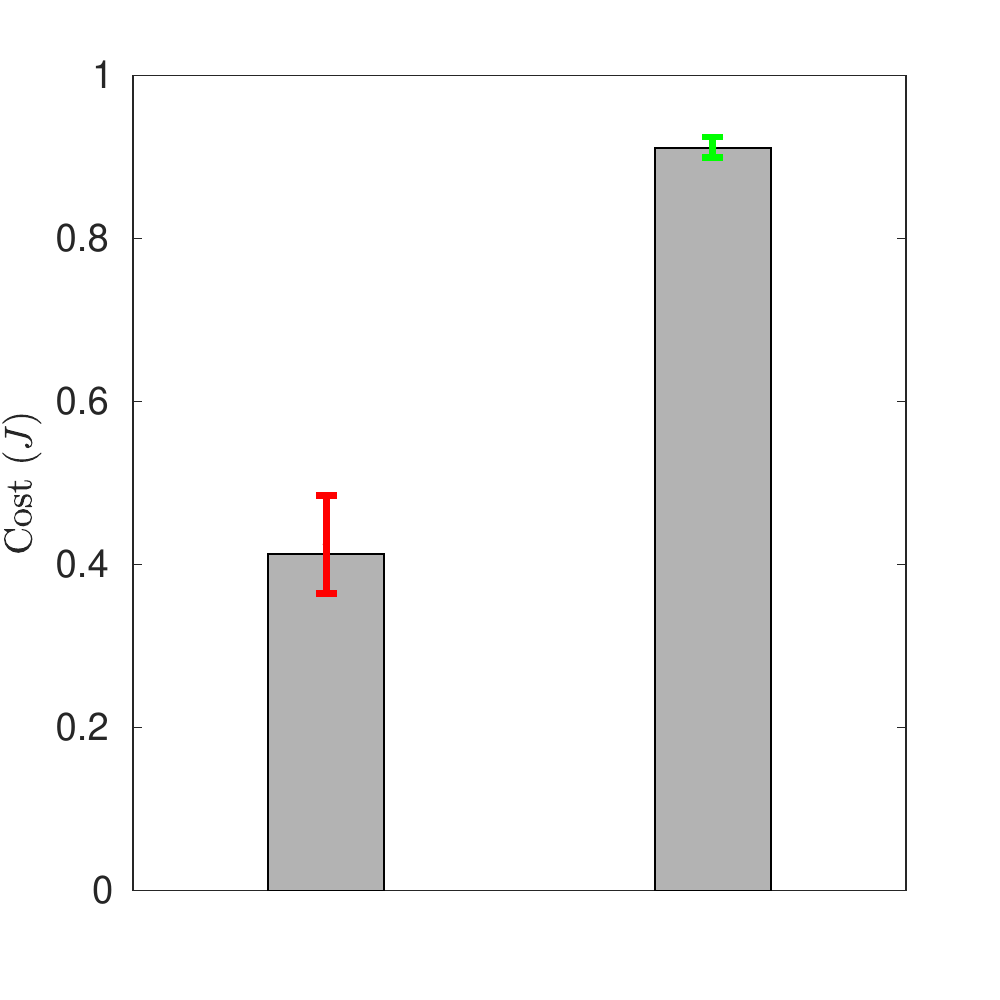}}  
	\caption{Results obtained for an LQR problem with a stable $A$ matrix being uncertain}
	\label{fig:LS_a}
\end{figure}

The results for the case where $a$ is the uncertain parameter with its nominal value as $a_0 = -1$ (stable), and $b = 1$ are shown in Figure \ref{fig:LS_a}.
Since $a$ is the source of uncertainty, by switching on the desensitization ($Q = 1,000$), it can be observed from Figure~\ref{fig:LS_a}(a) that the state approaches zero faster compared to the non-desensitized solution.
Consequently, from the Monte-Carlo simulations ($a \in [0.8a_0, 1.2a_0]$), it can be observed that the variations in the optimal trajectory (Figure \ref{fig:LS_a}(c)), and the cost (Figure \ref{fig:LS_a}(d)) are significantly lower for the desensitized solution, though the cost for the same is higher which is a trade-off. 
The error bars in Figure \ref{fig:LS_a}(d) represent the minimum and the maximum costs obtained form the Monte-Carlo results where the corresponding grey bars represent the nominal costs with $a=a_0$.

\begin{figure}[htb!]
	\centering
	\subfigure[Monte-Carlo simulations]{\includegraphics[width=0.4\textwidth]{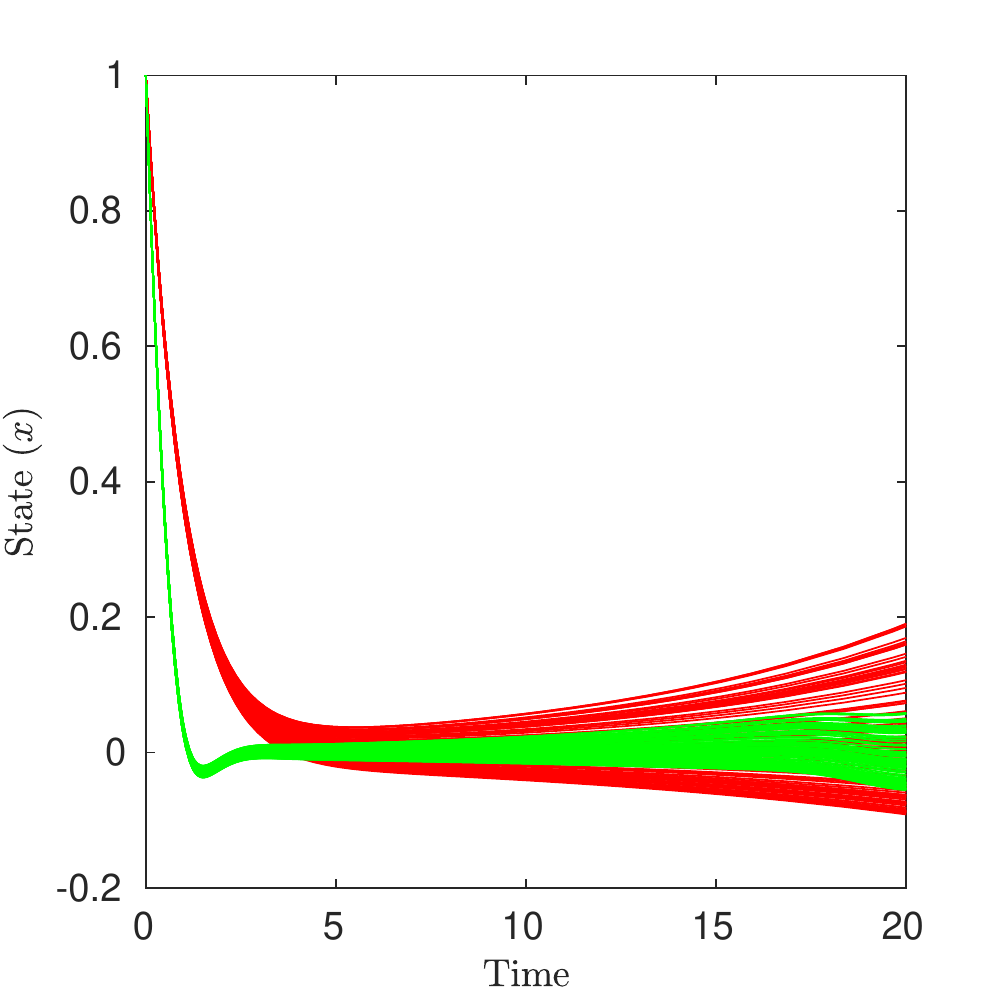}}
	\subfigure[Cost ($J$)]{\includegraphics[width=0.4\textwidth]{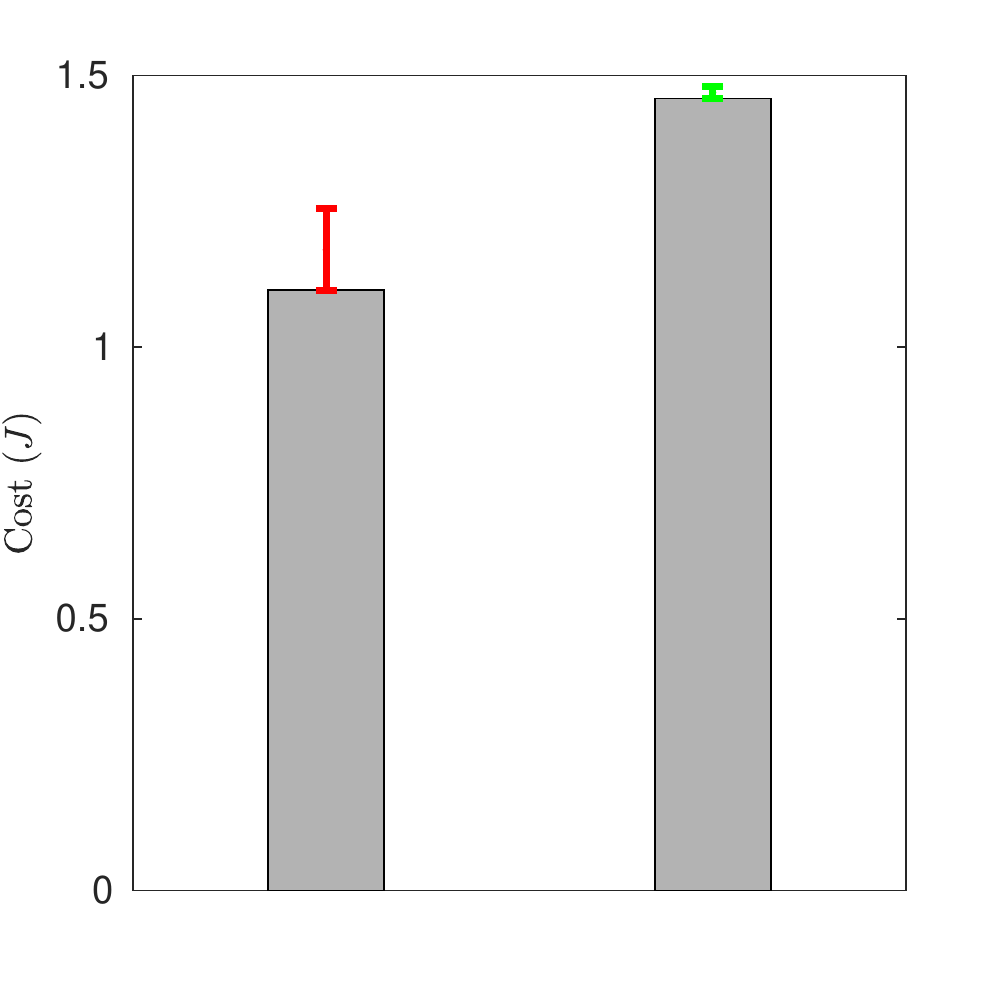}}   
	\caption{Results obtained for an LQR problem with an unstable $A$ matrix being uncertain (red: $Q = 0$, green: $Q = 1,000$).}
	\label{fig:LS_a1}
\end{figure}

\begin{figure}[htb!]
	\centering
	\subfigure[Monte-Carlo simulations]{\includegraphics[width=0.4\textwidth]{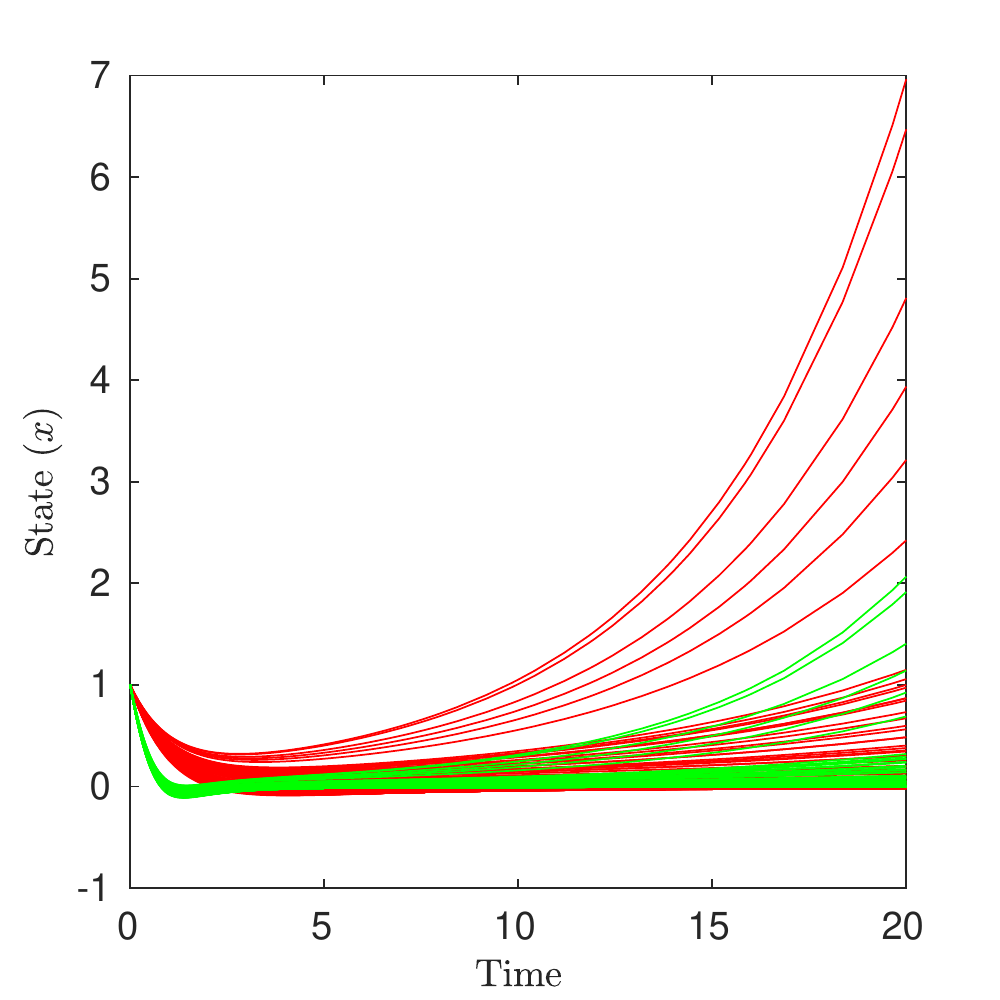}}
	\subfigure[Cost ($J$)]{\includegraphics[width=0.4\textwidth]{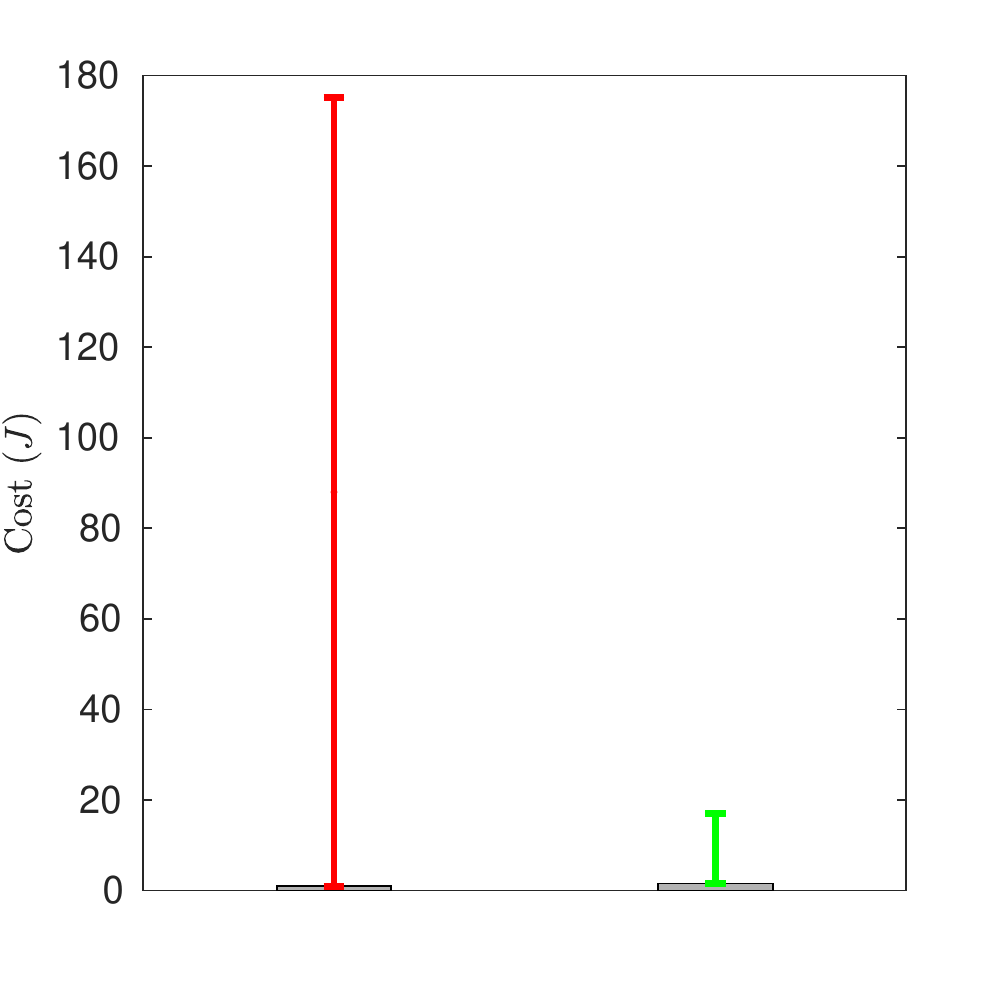}}   
	\caption{Results obtained for the LQR problem with $a = 0$ (red: $Q = 0$, green: $Q = 1,000$).}
	\label{fig:LS_a0}
\end{figure}

We also consider an unstable system with $a\in [0.8a_0,1.2a_0]$ being the uncertain parameter with $a_0 = 0.1$. 
Figure \ref{fig:LS_a1}(a) shows the Monte-Carlo results, and the behavior is similar to the one in the previous example except that the trajectories are now divergent because the system matrix is unstable.

A more interesting case is a marginally stable system with $a_0 = 0$, and $a\in [-0.2,0.2]$.
The corresponding results can be found in Figure \ref{fig:LS_a0}.
In the previous cases, although a parametric variation in $a$ is studied, such variations did not change the stability of the system, i.e., if the nominal system is stable, then the system with parametric variation is stable as well. 
Since $a$ can be both stable and unstable, the optimal control obtained for the nominal system without desensitization will be 
less effective combating the instabilities compared to the desensitized solution, as can be seen from the dispersion in trajectories (and costs) in Figure~\ref{fig:LS_a0}. 


\section{Discussion: Relation between the Sensitivity Matrix and co-states}
\label{sec:relation}

In this section we address the relationship between the sensitivity matrix defined in \eqref{Sey:SF} and the co-states $\lambda$.
Let us note that,
\begin{align}
\lambda^\top(t) = \frac{\partial C}{\partial x} (x(t),u,t),
\end{align}
which can be expressed as
\begin{align}
\lambda^\top(t) &= \frac{\partial C(x(t),u,t)}{\partial x_0}\left[ \frac{\partial x(t)}{\partial x_0}\right]^{-1} \nonumber  \\
&= \frac{\partial C(x(t),u,t)}{\partial x_0} S(t|t_0,x_0)^{-1}, \label{eq:co-state_sm_rel}
\end{align}
where $S(t|t_0,x_0) = \dfrac{\partial x(t)}{\partial x_0}$ is the sensitivity of the state at time $t$ along the trajectory with respect to variation in its initial condition $x_0$ \cite{seywald1996desensitized}.
Note that the dependency of $x(t)$, $t \geq t_0$ on $t_0$ and $x_0$ (initial conditions) is implicit.
The relationship between the co-state and the sensitivity matrix in (\ref{eq:co-state_sm_rel}) can be generalized to obtain the sensitivity of the solution with regards to the state at any other time $t'$ as

\begin{align}
\lambda^\top(t) &= \frac{\partial C(x(t),u,t)}{\partial x(t')}\frac{\partial x(t')}{\partial x(t)}, \\
&= \frac{\partial C(x(t),u,t)}{\partial x(t')} S(t|t',x(t'))^{-1}~\forall\ t,t' \in [t_0,t_f].
\end{align}
Therefore,
\begin{align}
\lambda(t) &= S(t|t',x(t'))^{-\top} \left[ \frac{\partial C(x(t),u,t)}{\partial x(t')} \right]^\top,\\
\left[ \frac{\partial C(x(t),u,t)}{\partial x(t')} \right]^\top &= S(t|t',x(t'))^{\top} \lambda(t).
\end{align}

From the above expressions, we observe that the sensitivity matrix $S(t|t',x(t'))^{\top}$ is essentially the transition matrix between the co-states $\lambda(t)$ and the partial of the cost-to-go function at time $t$ with respect to the state at time $t'$, i.e., the sensitivity of the cost-to-go function at time $t$ with respect to the state at time $t' < t$.


\section{Conclusion} \label{sec:conclude}

We attempt to exploit the co-states to obtain trajectories that are \emph{less sensitive} to parametric variations.
It is established that the co-states, defined by the Hamiltonian and the adjoint equations, capture the sensitivity of a cost-to-go function 
for \textit{any} arbitrary control law. 
This has led to the idea of inserting the co-states into the cost function and then studying its implications in the context of DOC.
In particular, this approach has been used to solve the problem of cost desensitization with respect to variation in system parameters.
The results suggest that variations to parametric uncertainty and optimality can be balanced using this approach through the choice of an appropriate weighting parameter. 
The numerical simulations give promising results that validate the theory. 
The proposed approach can be used to look at some other interesting problems, especially related to robust optimal control. 
It is also of interest to closely study its connections with covariance-steering for stochastic dynamical systems~\cite{okamoto2018Optimal}.


\section*{Acknowledgment}

This work has been supported by NSF award CMMI-1662542.


\section*{Appendix}

\subsection{Proof of Theorem \ref{thm:main}} \label{Ap:proof}


For a fix control $\bar{u}$, the cost-to-go  from any state $x$ at time $t$  is
\begin{align*}
    C(x,\bar u,t)=&\phi(x(t_f),t_f)+\int_{t}^{t_f}L(x(s), \bar{u}(s),s)\mathrm ds
\end{align*}
where $x(t)=x$.

Let the perturbed state at time $t$  be represented by $x(t,\alpha)=x(t)+\alpha \delta x (t)$ where $\alpha \in [0,\alpha_0)$ for some $\alpha_0>0$, and $\delta x (t)\in \mathbb R^n$.
With this perturbation the new cost-to-go is 
\begin{align*}
C(x+\alpha \delta x,\bar{u},t)=  &  \phi(x(t_f,\alpha),t_f)+\\
&~~\int_{t}^{t_f}L(x(s,\alpha),\bar{u}(s),s)\mathrm ds
\end{align*}
where $x(s,\alpha)$ denotes the perturbed state at time $s\ge t$.
By denoting 
\begin{align*}
    x(s,\alpha)=x(s,0)+\alpha \mathrm \delta x(s),
\end{align*}
for all $s\ge t$, we obtain
\begin{align*}
    \delta \dot{x}(s)=f_x(x(s,0),\bar{u}(s),s)\mathrm \delta x(s)+O(\alpha), 
\end{align*}
where $O(\alpha)$ is such that $\lim_{\alpha\to 0}O(\alpha)=0$. 
Consequently, $\delta x(s)=\Gamma(s,t)\delta x(t) +O(\alpha)$
where $\Gamma(s,t)$ is the state transition matrix corresponding to the matrix $f_x(x(s,0),u(s),s)$.

Therefore,
\begin{align*}
    C(x+&\alpha \delta x,\bar{u},t)-C(x,\bar{u},t)=
    \alpha\phi_x(   x(t_f,0), t_f) \Gamma(t_f,t)\delta x +\alpha\left[\int_{t}^{t_f}L_x(x(s,0),u(s),s)\Gamma(s,t)\mathrm ds\right]\delta x
    +O(\alpha^2)
\end{align*}
and thus,
\begin{align*}
    \lim_{\alpha\to 0^+}\frac{ C(x+\alpha \delta x,\bar{u},t)-C(x,\bar{u},t)}{\alpha} &= \Big[\phi_x(x(t_f,0),t_f)\Gamma(t_f,t)+ \int_{t}^{t_f}\,L_x(x(s,0),u(s),s)\Gamma(s,t)\mathrm ds\Big]\delta x,\\
   \implies C_x(x,\bar u,t)& =\phi_x(x(t_f,0),t_f)\Gamma(t_f,t) +\int_{t}^{t_f}L_x(x(s,0),\bar{u}(s),s)\Gamma(s,t)\mathrm ds.
\end{align*}
At this point, if we denote
\begin{align*}
    \lambda^{\top}\triangleq C_x(x,\bar{u},t).
\end{align*}
We then have
\begin{align*}
    \dot{\lambda}^{\top}=-\lambda^{\top}f_x-L_x,
\end{align*}
since $\dot \Gamma(s,t)=-\Gamma(s,t)f_x(x(t,0),u(t),t)$.
Furthermore,
$\lambda(s)$ satisfies the terminal condition $\lambda(t_f)=\phi_x(x(t_f),t_f)$.
Thus, if we define the Hamiltonian as $H=L+\lambda^{\sf T}f$, it follows that 
\begin{align*}
    \dot{\lambda}^{\top}=-\frac{\partial H}{\partial x}
\end{align*}
and this $\lambda$ represents the first order variation in cost-to-go with boundary condition
\begin{align*}
    \lambda(t_f)=\phi_x(x(t_f),t_f).
\end{align*}
The result follows.
\qed


\end{document}